\documentclass{article}

\usepackage{amsfonts}
\usepackage{amsmath}
\usepackage[latin1]{inputenc}
\usepackage{hyperref}
\title{Where do odd perfect numbers live?}
\author{Aldi Nestor de Souza\\ \url{aldinestor@ufmt.br}}

\begin{document}
\maketitle
\section{Abstract}
The existence of a perfect odd number is an old open problem of number theory. An Euler's theorem states that if an odd integer $ n $ is perfect,  then $ n $ is written as $ n = p ^ rm ^ 2 $, where $ r, m $ are odd numbers, $ p $ is a prime number of the form $ 4 k + 1 $ and $ (p, m) = 1 $, where $ (x, y) $ denotes the greatest common divisor of $ x $ and $ y $. In this article we show that the exponent $ r $, of $ p $, in this equation,  is necessarily equal to 1. That is, if $ n $ is an odd perfect number, then $ n $ is written as $ n = pm ^ 2. $
\section{Introduction}
A positive integer $ n$ is said perfect when $ \sigma (n) = 2n $, where $ \sigma (n) $ denotes the sum of the positive divisors of $ n $. For example, 6, 28, and 496 are perfect. All known perfect numbers are even and an Euclide's theorem characterizes such numbers. The theorem reads as follows: an even number $n$ is perfect if, and only if, $ n = 2^{p-1} (2^p-1) $, where $ 2 ^ p-1 $ is a Mersenne's prime number. It is not known if there are infinite even perfect numbers. The existence of odd perfect numbers, on the other hand, remains an open problem. There is a vast literature, replete with many curious results on the subject. For example, Steuerwald (1937) has shown that numbers of the form $ n = p ^ kp_1 ^ 2p_2 ^ 2 \cdots p_k ^ 2 $, where $ p, p_1, \cdots, p_k $ are prime numbers and $ k $ is an odd number,  are not perfect; according to [1], odd perfect numbers have at least nine prime factors; according to $ [2] $ odd perfect numbers are greater than $ 10 ^ {1500};$ according to $ [3], $ every odd perfect number has a prime factor that exceeds $ 10 ^ 6 $; according to $ [4] $, any odd perfect number has the form $ 12k + 1 $ or $ 36m + 9 $. In many of the results on odd perfect numbers, the equation $ n = p ^ km ^ 2 $, given by Euler's theorem, plays a fundamental role and it is with respect to this equation, according to the following theorem, the main result of this article.

\section{Development}

\subsection{Theorem} 
If $n$ is an odd perfect number, then $ n $ is written as $$ n = pm^2. $$

\textbf{Demonstration}\\

From Euler's theorem, see [5], if n is an odd perfect number, then $ n = p ^ rm ^ 2 $, where $ m, r $ are odd numbers, $ p $ is a prime number of the form $ p = 4k + 1 $ and $ (p, m) = 1 $, where $ (x, y) $ denotes the maximum common divisor of $ x $ and $y$. Thus we have $\sigma(p^rm^2)=2p^rm^2.$ That is, 
\begin{equation}(1+p+\cdots+p^r)\sigma(m^2)=2p^rm^2.\end{equation} Thus $p$ is a root of the equation \begin{equation}
(a-b)p^r+ ap^{r-1}+\cdots+ap+a=0, 
\end{equation} where $a=\sigma(m^2)$ and $b=2m^2.$ Dividing the two sides of the equation $(2)$, by $ a $, we have that $ p $ is a root of \begin{equation} (1 -\frac {b}{a}) x ^ r + x ^ {r-1} +\cdots + x + 1 = 0. 
\end {equation}
\newpage 
Note that, since $ m ^ 2 $ is not perfect, $ a \neq b $. Note also that if $ \sigma (m ^ 2)> 2m ^ 2 $, that is, if $ a> b $, then $ p $ will not be a root of the equation $ (3) $ since all coefficients on the left side of $(3)$ are positive. Thus we have $ m ^ 2 <\sigma (m ^ 2) <2m^2 $, that is $ 1 <\frac {2m ^ 2} {\sigma (m ^ 2)} <2. $ Therefore, $1<\frac{b}{a} <2.$ Now let $ d =(a, b), a_1 = \frac {a} {d}, b_1 = \frac {b} {d}. $ It follows that $ p $ is a root of the equation 
\begin{equation}(1-\frac{b_1}{a_1})x^r + x^{r-1} + \cdots +x + 1 = 0. \end{equation}

Note that, since $(a_1, b_1) = 1$, if an integer $s$ is a root of the equation $ (4) $, then $s = ta_1,$ for some $t\in \mathbb Z$. Supposing, by absurdity, that, $ r> 1 $, and assuming that $ s = ta_1 $ is a root of the equation $ (4) $, it follows that $$ (1- {b_1}) t ^ ra_1^{r-1} + t^{r-1}a_1 ^{r-1} +\cdots + ta_1 + 1 = 0. $$
It follows that $ a_1> 1 $ is a divisor of 1, which is an absurdity, and therefore no integer, in particular no prime number, is a root of $ (4) $ \rule{0.2cm}{0.2cm}\\

Due the above theorem, if an odd number $ n $ is perfect, then there are an odd number $ m $ and a prime number $ p $ such that $ n = pm ^ 2 $. As $(p + 1)\sigma(m^2) = 2pm^2 $, it follows that $p $ is a root of the linear equation $$ (\sigma(m^2) -2m ^ 2)x + \sigma(m^2)= 0.$$ That is, 
$$p=\frac{\sigma(m ^ 2)}{2m^2- \sigma(m ^ 2)}.$$

\end{document}